\documentclass[12pt]{article} \def\C{\mathbf{C}}
\def\R{\mathbf{R}}
\def\id{\mathrm{id}}
\def\bC{\mathbf{\overline{C}}}

\def\N{\mathbf{N}}
\def\Rea{\mathrm{Re\ }}
\def\Ima{\mathrm{Im\ }}
\begin{document}
\title{Invariant curves and semiconjugacies of rational functions}
\author{Alexandre Eremenko\thanks{Supported by NSF grant DMS-1067886.}}
\maketitle
\begin{abstract} Jordan analytic curves which are invariant under
rational functions are studied.

Keywords: rational functions, iteration, functional equations.

MSC 30D05, 37F10.
\end{abstract}

In the last paragraph of his memoir \cite{F}
on iteration of rational functions Fatou writes:

{\em ``Il nous resterait \`a \'etudier les courbes
analytiques invariantes par une transformation rationnelle
et dont l'\'etude est intimement li\'ee
\'a celle des fonctions \'etudi\'ees dans ce Chapitre.
Nous esp\'erons y revenir bient\^ot.}\footnote{It remains to
study analytic curves invariant under rational transformations, which are
intimately connected with the functions studied in this chapter.
We hope to return to this soon.}

As far as I know, Fatou never returned to this
question
in his published work. Neither I know of any systematic study of
the question after Fatou. If a Jordan analytic invariant curve is
the boundary of a domain of attraction, then it must be a circle, \cite{F}.

Which Jordan analytic curves  $\gamma$
in the Riemann sphere can be invariant under rational
functions? Of course, $\gamma$ can be a circle,
and it is easy to describe all rational
functions which leave a given
circle invariant: such functions must commute with the reflection
in this circle. 

Other examples are obtained as level lines of the linearizing functions
in Siegel discs or Hermann rings.
It is not known whether a Jordan analytic invariant curve, different
from a circle, and which is mapped onto itself homeomorphically
can intersect the Julia set.

If $f$ is a polynomial (or an entire function)
then the only possible Jordan invariant
analytic curves in $\C$ are either circles or preimages of circles
under linearizing functions of Siegel discs \cite{A}.

Przytycki asked whether a rational function $f$ can have a repeller,
which is a Jordan analytic curve different from a circle. A repeller $\gamma$
is a compact set which has a neighborhood $V$, such that for $z\in V$,
$f^n(z)\in V$ if and only if $z\in\gamma$.

The principal result of this paper implies that such Jordan analytic
repellers must be algebraic curves, except for the Latt\'es examples.
However, the only example of
a rational function having 
such a repeller, different from a circle, that I could produce,
is a Latt\'es example.
\vspace{.1in}

\noindent
{\bf Theorem 1.} {\em Let $f$ be a rational function and $\gamma$
a Jordan analytic invariant curve such that $f\vert_\gamma$ is
not a homeomorphism and there is a repelling periodic point
$a\in\gamma$ with multiplier $\lambda,\; |\lambda|>1$. Assume in addition
one of the following:

a) 
$\gamma$ does not contain critical or neutral rational
fixed points of $f$, or

b) $\gamma\subset J(f)$.

Then $\lambda$ is real, and there exist a non-hyperbolic
Riemann surface $S$
with an anti-conformal involution $s$ and an
endomorphism $g:S\to S$, and a holomorphic map
$h:S\to\bC$ such that $g\circ s=s\circ g$,
$h(X)\subset\gamma$ where $X\subset S$ is the set of points fixed by
the involution $s$, 
and the following semi-conjugacy relation holds:
\begin{equation}\label{semi}
h\circ g=f^n\circ h
\end{equation}
with some integer $n\geq 0$.
}
\vspace{.1in}

\noindent
{\bf Corollary.} {\em Assumptions of Theorem 1, imply that the curve
$\gamma$ is algebraic, or $f$ is a Latt\'es example.}
\vspace{.1in}

If $\gamma$ is a repeller, then all assumptions of Theorem 1,
including both a) and b) are satisfied.
\vspace{.1in}

{\em Proof of Theorem 1}. Replacing $f$ by some iterate, we reduce to the
case that the point $a$ is fixed. Let $F$ be the Poincar\'e
function associated with the
fixed point $a$. This means that
\begin{equation}\label{poin}
F(\lambda z)=f\circ F(z),\quad F(0)=a,\quad F'(0)\neq 0.
\end{equation}
Such function exists for every repelling fixed point and it is
meromorphic in the plane $\C$.
Let $\Gamma$ be the component of $F^{-1}(\gamma)$
that contains $0$. Then the intersection of $\Gamma$
with a neighborhood of $0$ is a smooth curve invariant under
the map $z\mapsto\lambda^{-1}z$, which implies that $\lambda$ is real.
Replacing $f$ by the second iterate we achieve that $\lambda>1$.
Then it follows that $\Gamma$ contains a straight line,
and without loss of generality we may assume that this is the real line.
We have $F(\R)\subset\gamma$. 

Consider the restriction mappings $\R \to \gamma$ of
all Poincar\'e maps $F_j$ for all repelling
fixed points $a_j$ on $\gamma$. There are finitely many of them,
as the number of fixed points is finite.
We claim that under our assumptions
at least one
of these maps is not injective.  

Indeed, suppose that some $F_j$ is injective on the real line.
Then the image $F_j(\R)$ is a simple arc $\gamma_j\subset\gamma$,
and it is easy to see that the endpoints of this arc must be 
attracting fixed points with real multipliers.

Indeed these endpoints are fixed because $\gamma_j$ is invariant
and the iterates of $f$ converge to these endpoints.
This convergence implies that these endpoints cannot be repelling.

Suppose that the assumption b) in Theorem 1 holds.
As $\gamma_j$ is a subset of a Jordan analytic curve, its endpoints
cannot be neutral irrational points. A local description of dynamics
near a neutral rational fixed point shows that if such a point $w$ is an
endpoint of a smooth invariant curve on which the iterates converge
to $w$, then this curve must intersect the set of normality.
This contradiction
completes
the proof of non-injectivity of $F_j\vert_\R$ under the assumption b).

We continue the non-injectivity proof
under the assumption a). 
The endpoints of $\gamma_j$ can coincide
in which case we obtain that $f\vert_\gamma$ is a homeomorphism.
Now let $b$ be an endpoint of $\gamma_j$ which does not coincide
with the other endpoint. Assumption a) implies that the endpoints
are attracting. From the local description of the dynamics near an
attracting fixed point $b$
we conclude that there is a simple invariant arc $\delta_0$ disjoint
from $\gamma_j$ such that $\gamma_j\cup\delta_0\cup\{ b\}$ is
an arc of $\gamma$. For $k\geq 1$,
let $\delta_k$ be the component of $f^{-1}(\delta_{k-1})$
which contains $\delta_{k-1}$. Then
$$\Delta=\cup_{k=1}^\infty\delta_k$$
is a simple invariant curve, $\Delta\subset\gamma$,
one endpoint of $\Delta$ is $b$ and the other endpoint is a repelling
fixed point $a_k$ on $\gamma$ distinct from $a_j$. This means 
that the curves $\gamma_j$ and $\gamma_k$ have a common endpoint $b$,
and their union is an analytic curve near this common point. 

We conclude that our curves $\gamma_j$ cover $\gamma$, and $f\vert_\gamma$
is a homeomorphism. This contradicts the assumptions and proves the
claim.

From now on we assume that some Poincar\'e function $F$
has non-injective restriction on the real line.

We recall a result of \cite{ER}.
Let $F:\C\to S$ be a non-constant holomorphic map from the complex plane
to a Riemann surface $S$.
Consider the following equivalence relation in $\C$:
$x\sim y$ iff $F(x)=F(y)$. Let $\Gamma_F\subset\C^2$ be the graph of this
equivalence relation. It is easy to see that
$\Gamma_F$ is an analytic variety of pure (complex) dimension $1$
(has no isolated points). 

Now suppose that a analytic variety $\Gamma\subset\C^2$ of pure dimension $1$
which is a graph of an equivalence relation is given.
Then there exists a holomorphic map $F:\C\to S$,
where $S$ is a non-hyperbolic Riemann surface such that 
$F(x)=F(y)$ iff $(x,y)\in \Gamma$. This map $F$ is defined
by $\Gamma$ uniquely up
to a composition with an automorphism of $S$.

Now we characterize the graphs $\Gamma$ of equivalence relations
corresponding to Poincar\'e functions.
\vspace{.1in}

\noindent
{\bf Lemma 1.} {\em A holomorphic map $F:\C\to S$ to a non-hyperbolic Riemann
surface $S$ is a 
Poincar\'e function
of an endomorphism of $S$ if and only if
$\Gamma_F$ is invariant under a map 
\begin{equation}
\label{lambda}
(x,y)\mapsto(\lambda x,\lambda y)
\end{equation}
for
some $\lambda\in\C, |\lambda|>1$.}
\vspace{.1in}

{\em Proof.} That the graph $\Gamma_F$ corresponding to a function $F$
satisfying (\ref{poin}) is invariant with respect to this map is clear.

To prove the converse statement, we recall the result from \cite{ER}
that the existence of a decomposition $F=f\circ G,$
where $f$ and $G$ are maps between
non-hyperbolic Riemann surfaces is equivalent to the inclusion
$\Gamma_G\subset\Gamma_F$. Let $\Gamma_1=\lambda\Gamma_F\subset\Gamma_F.$
Then $\Gamma_1=\Gamma_G$ where $G(\lambda z)=F(z)$.
On the other hand, by the result of \cite{ER} just cited,
$F=f\circ G$ where $f:S\to S$ is an endomorphism.
So $F(\lambda z)=f\circ G(\lambda z)=f\circ F(z)$, that is $F$ satisfies
a Poincar\'e equation. Putting $z=0$ we obtain that $F(0)$ is a fixed
point of $f$, and the multiplier of this fixed point is $\lambda$ by
the chain rule.

This completes the proof of the lemma.
\vspace{.1in}

A holomorphic map $F:\C\mapsto S$ will be called real 
if there exists an anti-conformal involution $s:S\to S$
such that $F(\overline{z})=s\circ F(z)$.
\vspace{.1in}

\noindent
{\bf Lemma 2.} {\em  A holomorphic map
$F:\C\to S$ is real if and only if $\Gamma_F$
is invariant under the map 
\begin{equation}\label{conj}
(x,y)\mapsto(\overline{x},\overline{y}).
\end{equation}}

This is clear.

Now we return to the proof of Theorem 1. As $F$ is not injective on the 
real line, there is a real analytic germ $\phi\neq\id$
such that $F\circ\phi=F$.
This implies that some part $\Gamma_1\neq\{(x,x):x\in\C\}$
of the graph $\Gamma_F$ is parametrized
by $(x,\phi(x))$, so this part is invariant under the map (\ref{conj}).

Now let $\Gamma_2\subset\Gamma_F$ be the smallest analytic variety
of pure dimension $1$ which contains $\Gamma_1$, which is a graph of
an equivalence relation and which is invariant under the maps 
(\ref{conj}) and (\ref{lambda}).
\vspace{.1in}

{\em Proof of existence of $\Gamma_2$}. Let $E\subset\C^2$ be an
arbitrary set containing the diagonal $D=\{(x,x):x\in\C\}$.
Consider three operations on such sets:
$$E\mapsto\Lambda E=\{\lambda^n(x,y):(x,y)\in E,\; n\in\N\},$$
$$E\mapsto SE=\{(x,y):(y,x)\in E\},$$
$$E\mapsto TE=\{(x,z):\exists y, (x,y)\in E, (y,z)\in E\}.$$
If $E$ does not contain isolated points then $\Lambda E, SE,TE$ do not contain
isolated points. Moreover, if $E$ is symmetric with respect to the operation
(\ref{conj}) then each $\Lambda E, ES, TE$ is also symmetric with respect
to this operation. Now we apply all finite sequences of
operations $\Lambda, S, T$ to 
to $\Gamma_1\cup D$, and take as $\Gamma_2$ the union of those
irreducible components of $\Gamma_F$ that have non-isolated
intersection with the resulting
sets. Then $\Gamma_2$ is the minimal analytic set of pure dimension 1,
which is a graph of an equivalence relation, invariant with respect to
(\ref{lambda}), and by the previous remark, it is symmetric
with respect to (\ref{conj}). This completes the construction of
the set $\Gamma_2$.

If $\Gamma_2=\Gamma_F$ then $F$ is real and $\gamma$ is a circle.
If $\Gamma_2\neq\Gamma_F$ then there exists a factorization
\begin{equation}\label{fact}
F=h\circ G,
\end{equation}
where $G$ is a Poincar\'e function of some endomorphism $g:S\to S$
of a non-hyperbolic surface $S$,
that is 
\begin{equation}\label{G}
G(\lambda z)=g\circ G(z),
\end{equation}
$\Gamma_2=\Gamma_G$,
and $h: S\to\bC$ is a holomorphic map. Moreover, $G$ is real.
Combining (\ref{fact}) and (\ref{G}) we obtain
$$f\circ h\circ G(z)=f\circ F(z)=F(\lambda z)=h\circ G(\lambda z)=
h\circ g\circ G(z),$$
and this implies (\ref{semi}).
This completes the proof of Theorem 1.

\vspace{.1in}

To prove the Corollary,
we discuss the functional equation (\ref{semi}).
First of all, $S$ can be a torus, $\C^*, \C$ or $\bC$.
If $S$ is a torus, then $f$ is a Latt\'es example.
If $S$ is $\C$ or $\C^*$, and $h$ has an essential singularity
then $f$ also must be a Latt\'es example. This was proved in \cite{E}.
See also \cite{BE} for another proof.
Otherwise $S=\bC$, and
thus both $g$ and $h$ are rational. This proves the corollary.
\vspace{.1in}

Latt\'es examples indeed have Jordan analytic invariant curves
which are not circles. These curves can be algebraic or transcendental.
\vspace{.1in}

\noindent
{\em Example 1.} Let $\wp$ be the Weierstrass function with periods 
$2\omega_1$ and $2i\omega_2$ where we assume that
$\omega_1$ and $\omega_2$ are real.
Consider the line $L=\{ x+2i\omega_2/3: x\in\R\}$ and let $\gamma=\wp(L)$.
The simplest Latt\'es function $f$ corresponding to $\wp$ satisfies
\begin{equation}\label{lat}
\wp(2z)=f\circ\wp(z),
\end{equation}
and we see that $\gamma$ is invariant under $f$, because $2L\equiv-L$
modulo periods, and $\wp$ is even.
It is easy to see that $\gamma$ is a Jordan analytic curve
which is not a circle, and $f:\gamma\to\gamma$ is a two-sheeted covering
map. 

Let us show that $\gamma$ is algebraic. Let $s(z)$ be the reflection
in the line $L$. Define 
$$X(z)=(\wp(z)+\overline{\wp(s(z))})/2,$$
and
$$Y(z)=(\wp(z)-\overline{\wp(s(z))})/(2i).$$
Then $X$ and $Y$ are elliptic functions with the same period lattice as $\wp$.
So they are related by an algebraic equation
$$F(X,Y)=0.$$
When $z\in L$, we have $s(z)=z$, so $X(z)=\Rea\wp(z)$ and $Y(z)=\Ima\wp(z)$.
So points $x+iy$ on our curve $\gamma=\wp(L)$ satisfy
the algebraic equation $F(x,y)=0$.

\vspace{.1in}
\noindent
{\em Example 2.} Let $\wp$ be the Weierstrass
function with primitive periods
$(1,\tau)$, where $\tau=p+i,$ where $p$ is real and irrational.
Let $\gamma=\wp(L)$, where $L$ is the same as Example 1 above,
and $f$ is a Latt\'es function as in
(\ref{lat}).
It is easy to see that $\gamma$ is a Jordan analytic curve invariant 
under $f$, which is not a circle.

Let us show that $\gamma$ is transcendental.
Consider the function $g(z)=\overline{\wp(s(z))}$,
where $s$ is the reflection with respect to $L$. It is an
elliptic function with periods $1$ and $\overline{\tau}=p-i$.
Evidently, the period $p-i$ is not a rational combination
of $1$ and $p+i$.
So $g$ is an elliptic function whose lattice is not commensurable with
the lattice of $\wp$. Suppose now that $\gamma$ is algebraic and
let $F(x,y)=0$ be the equation of $\gamma$. Then 
$$F((\wp+g)/2,(\wp-g)/(2i))=0$$
holds on the real line and thus everywhere.
We conclude that $\wp$ and $g$ are
algebraically dependent, and this is a contradiction, because algebraically dependent
elliptic functions must have commensurable lattices. This proves that
$\gamma$ is transcendental.

\vspace{.1in}

Equation (\ref{semi}) in rational functions was recently studied
in \cite{I} and  \cite{Pako}. 
It is closely related to factorization theory of rational
functions with respect to composition which is due to J. Ritt
\cite{R1,R2,R3}.
Pakovich \cite{Pako} obtained a complete classification of
triples of rational functions that satisfy (\ref{semi}).
However, it is hard to obtain from this classification a description
of triples that satisfy other conditions of Theorem 1.

The simplest solutions of (\ref{semi}) with $n=1$ can be constructed
as follows \cite{I}:

Let $u$ and $v$ be two rational functions. Set
\begin{equation}\label{triple}
f=u\circ v, \quad g=v\circ u \quad \mbox{and}\quad h=u.
\end{equation}
Then $h\circ g=f\circ h$, so (\ref{semi}) is satisfied.

Another class of examples is obtained by taking an arbitrary rational
function $w$ and setting $f(z)=z^mw^n(z)\; g(z)= z^mw(z^n)$ and
$h(z)=z^n$.

These examples do not exhaust all possibilities which are listed
in \cite{Pako}.

It is interesting to know whether (\ref{triple}) leads to Jordan analytic
invariant curves which are not circles.
Such examples will occur if 
\begin{equation}\label{factor}
v\circ u
\end{equation}
is a real rational function
but $u$ is not real and maps the circle $\R\cup\{\infty\}$
onto a Jordan analytic curve $\gamma$ which is not a circle.
Then $\gamma$ will be an
invariant curve for $u\circ v$.
\vspace{.1in}

\noindent
{\em Example 3.} (F. Pakovich). Let $P_n(z)=z^n$ and $J(z)=(z+z^{-1})/2.$
It is well known that
$$R:=J\circ P_n=T_n\circ J,$$
where $T_n$ are real polynomials (they are Chebyshev polynomials normalized
so that the leading coefficient is $2^{n-1}$).
Now, if $\epsilon=\exp(2\pi i/n)$, then the first factorization implies
that $R(\epsilon z)=R(z)$.
So $R(z)= T_n\circ J(\epsilon z).$
Let $u(z)= J(\epsilon z)$.
Then $u(\R)$ is a hyperbola $\gamma$, it is mapped by $T_n$
on the real line,
and this hyperbola is invariant
under the map $f=u\circ T_n$.

However a hyperbola, when considered on the Riemann sphere,
is not a Jordan curve. Thus the only examples we have to illustrate
Theorem 1 are Examples 1 and 2\footnote{After this paper was accepted,
Peter M\"uller \cite{EM} found examples of rational functions $f_1$ and $f_2$,
such that $f_1\circ f_2$ is real but
$\gamma=f_2(\R)$ is a Jordan analytic curve
which is not a circle. This curve $\gamma$ is invariant under $f_2\circ f_1$,
thus giving a non-trivial illustration to Theorem 1, and an example
of a Jordan analytic repeller.
His functions $f_1$ and $f_2$
come from isogenies of elliptic curves. Actually in M\"uller's example,
both $f_1\circ f_2$ and $f_2\circ f_1$ are Latt\'es maps, but one
can consider an arbitrary real rational function $p$ and set $f_3=p\circ f_1$.
Then $g=f_3\circ f_2$, $f=f_2\circ f_2$ and $h=f_2$ will satisfy (\ref{semi})
and none of them will be a Latt\'es map.}.

The semiconjugacy equation (\ref{semi}) for rational functions
occurs also in \cite{J} and \cite{G} in different contexts.

The author thanks the referee for valuable remarks.

{\em Department of mathematics, Purdue University, West Lafayette, IN 27907}
\end{document}